\magnification\magstephalf

\overfullrule=0pt
\def\bE{\hbox{\bf E}}
\def\bR{\hbox{\bf R}}
\pageno=1

\centerline{\bf Remarks on Talagrand's deviation inequality for}
\centerline{\bf  Rademacher functions}\vskip6pt
\centerline{William B. Johnson*}
  \vfootnote*{Supported in part by NSF DMS-8703815.}

\centerline{Department of Mathematics}
\centerline{ Texas A\&M University}
\centerline{ College Station, TX
77843}
\centerline{and}
\centerline{Gideon Schechtman}
\centerline{Department of Mathematics}
\centerline{Texas A\&M University}
\centerline{Department of Theoretical Mathematics}
\centerline{The Weizmann Institute of Science}
\centerline{Rehovot, Israel}

\baselineskip=18pt

\bigskip

\noindent {\bf Introduction.}

Recently Talagrand [T] estimated the deviation of a function on $\{0,1\}^n$ from
its median in terms of the Lipschitz constant of a convex extension of $f$ to
$\ell ^n_2$; namely, he proved that

$$P(|f-M_f| > c) \le 4 e^{-t^2/4\sigma ^2}$$

\noindent where $\sigma$ is  the Lipschitz constant of the extension of $f$ and
$P$ is the natural probability on $\{0,1\}^n$.

Here we extend this inequality to more general product probability spaces; in
particular, we prove the same inequality for $\{0,1\}^n$ with the product
measure $((1-\eta)\delta _0 + \eta \delta _1)^n$. We believe this should be
useful in proofs involving random selections. As an illustration of possible
applications we give a simple proof (though not with the right dependence on
$\varepsilon$) of the Bourgain, Lindenstrauss, Milman result [BLM] that for
$1\le r < s \le 2$ and $\varepsilon >0$, every $n$-dimensional subspace of $L_s
\ (1+\varepsilon)$-embeds into $\ell ^N_r$ with $N = c(r,s,\varepsilon)n$.

\noindent {The main results.}

For $i=1,\ldots, n$ let $(X_i, \|\cdot \|_i)$, be normed spaces, let $\Omega _i$
be a finite subset of $X_i$ with diameter at most one and let $P_i$ be a
probability measure on $\Omega _i$.  Define

$$\eqalignno{X &= \bigg( \sum ^n_{i=1} \oplus X_i\bigg)_2\cr
\noalign{\hbox{and}}
\Omega &= \Omega _1 \times \Omega _2 \times \cdots \times \Omega _n \subset X}$$

\noindent and let

$$P = P^{(n)} = P_1 \times P_2 \times \cdots \times P_n$$

\noindent  be the product probability measure on $\Omega$. For a subset
$A\subseteq \Omega$ and $t\in \Omega$ let

$$\phi _A(t) = d(t, \ {\rm conv} \ A)$$

\noindent be the distance in $X$ from $t$ to the convex hull of the set $A$.

\proclaim Theorem 1. $\bE e^{{1\over 4} \phi ^2_A(t)} \le {1\over P(A)}$.

\noindent {\bf Remark 2:} \ Talagrand's theorem is the special case of Theorem 1
when each $\Omega _i$
consists of two points and $P_i$ gives weight ${1\over 2}$ to each of them.  In
the application below we use two point spaces for each $\Omega_i$, but $P_i$
does not assign the same mass to both points.\vskip12pt

\noindent {\bf Proof:} \ We repeat Talagrand's induction argument [T]; the
difference is only on the calculus level.  For $n=1$

$$\bE  e^{{1\over 4} \phi ^2_A(t)} \le P(A) + (1-P(A))e^{1\over 4} \le
{1\over P(A)},$$

\noindent as the maximal value of $r(r + (1-r)e^{1\over 4})$ for $0\le r
\le 1$ is 1.  Assume the theorem holds for $n$ and suppose that

$$A\subseteq \Omega _1\times \cdots \times \Omega _n \times \Omega _{n+1}.$$

\noindent Set, for $w\in \Omega _{n+1}$,

$$A_w = \{t\in \Omega _1 \times \cdots \times \Omega _n\colon \ (t,w) \in A\},$$

\noindent where for $t = (t_1,\ldots, t_n) \in \Omega _1 \times \cdots \times
\Omega _n$ and $w\in \Omega_{n+1}, (t,w)$ denotes $(t_1,\ldots, t_n,w)$.  Let $v
\in \Omega _{n+1}$ be such that

$$P^{(n)}(A_v) = \max_{w\in \Omega _{n+1}}P^{(n)}(A_w).$$

\noindent We shall use the following two inequalities:

$$\leqalignno{\phi ^2_A(t,v) &\le \phi ^2_{A_v}(t) \quad {\rm for \ all} \qquad
t\in \Omega _1\times \cdots \times \Omega _n&(1)\cr
\phi ^2_A(t,w) &\le \inf _{0\le \alpha \le 1}[\alpha \phi ^2_{A_w}(t) +
(1-\alpha) \phi ^2_{A_v}(t) + (1-\alpha)^2]&(2) \cr
&\hskip1in {\rm for \ all} \ t \in \Omega _1\times \cdots \times \Omega _n \
{\rm and \ for} \ w\ne v.}$$

\noindent Using (1) and (2), H\"older's inequality and the induction hypothesis
(in this order) we get

$$\leqalignno{\bE e^{{1\over 4}\phi ^2_A(s)} &= P_{n+1}\{v\}
\bE e^{{1\over 4}\phi ^2_{A_v}(t)} + \sum _{w\ne v} P_{n+1}\{w\} \inf _{0\le
\alpha \le 1} \bE \Big(e^{{1\over 4}\phi ^2_{A_w}(t)}\Big)^\alpha
\Big(e^{{1\over 4}\phi ^2_{A_v}(t)}\Big)^{1-\alpha} e^{{1\over
4}(1-\alpha)^2}\cr
&\le P_{n+1}\{v\} {1\over P^{(n)}(A_v)} + \sum _{w\ne v} P_{n+1}\{w\} \inf
_{0\le \alpha \le 1} \Big({1\over P^{(n)}(A_w)}\Big)^\alpha \Big({1\over
P^{(n)}(A_v)}\Big)^{1-\alpha} e^{{1\over 4}(1-\alpha)^2}&(3)\cr
&= {1\over P^{(n)}(A_v)} \Big[P_{n+1}(v) + \sum _{w\ne v} P_{n+1}(w) \inf _{0\le
\alpha \le 1} \Big( {P^{(n)}(A_v)\over P^{(n)}(A_w)}\Big)^\alpha e^{{1\over
4}(1-\alpha)^2}\Big].}$$

\noindent For $0\le \lambda \le 1$, let $\alpha (\lambda)$ be the point where
$\min\limits _{0\le \alpha \le 1} {1\over \lambda ^\alpha} e^{{1\over
4}(1-\alpha)^2}$ in attained; i.e.,

$$\alpha (\lambda) = \left\{\eqalign{&1+ 2 \log \lambda\cr &0}\right.
\eqalign{&{\rm if} \ 2 \log \lambda > -1\cr &{\rm otherwise}}$$

\noindent and set

$$g(\lambda) = {1\over \lambda ^{\alpha (\lambda)}} e^{{1\over
4}(1-\alpha(\lambda))^2} = \left\{ \eqalign{&e^{-\log \lambda - (\log \lambda
)^2}\cr &e^{1\over 4}}\right. \eqalign{&{\rm if} \ 2 \log \lambda > -1\cr &{\rm
otherwise.}}$$

\noindent From (3) we get

$$\leqalignno{\bE e^{{1\over 4} \phi ^2_A(s)} &\le {1\over P^{(n)}(A_v)}
\bigg[ P_{n+1}(v) + \sum _{w\ne v} P_{n+1}(w) g(\lambda _w)\bigg]&(4)\cr
\noalign{\hbox{where}}
\lambda _w &= {P^{(n)}(A_w)\over P^{(n)}(A_v)}.}$$

\noindent {\bf Claim:}	\ $g(\lambda) \le 2 -\lambda$ for $0\le \lambda \le
1$.\vskip9pt

  Using the claim we get from (4) that

$$\bE e^{{1\over 4}\phi ^2_A(s)} \le {1\over P^{(n)}(A_v)} (q +
(1-q)(2-t))\leqno (5)$$

\noindent where $q = P_{n+1}(v)$ and $t = {P(A) - P_{n+1}(v)P^n(A_v)\over
(1-P_{n+1}(v))P^n(A_v)}$ (note that $0\le q, t\le 1)$.	As

$${1\over P(A)} = {1\over P^{(n)}(A_v)} \ {1\over q+(1-q)t},\leqno (6)$$

\noindent it suffices to prove that

$$q + (1-q)(2-t) \le {1\over q+(1-q)t}\leqno (7)$$

\noindent for all $0 \le q, t\le 1$, which is easily checked.

The proof of the
claim is elementary: \ Let

$$f(\lambda) = g(\lambda) + \lambda -2.$$

\noindent Then $f(1) = f'(1) = 0$ and $f''(\lambda) \le 0, \ 0\le \lambda \le
1$,
which implies that $f(\lambda )\le 0$ for $0\le \lambda \le 1$.

\noindent {\bf Remark 3:} \ Let $2 < p < \infty$ and consider $\Omega$ as a
subset of $(\sum\limits ^n_{i=1} \oplus X_i)_p$. Set

$$\phi _{A,p}(t) = \inf \bigg\{ \bigg( \sum ^n_{i=1} \|t_i-s_i\|^p_i\bigg)
^{1/p}; s = (s_1,\ldots, s_n) \in \ {\rm conv} \ A\bigg\}.$$

\noindent Then, as pointed out by Talagrand, we also get

$$\bE e^{{1\over 4} \phi ^p_{A,p}(t)} \le {1\over P(A)}$$

\noindent because $\phi ^p_{A,p} (t) \le \phi ^2_A(t)$.

\proclaim Corollary 4. Let $2 \le p < \infty$ and let $f$ be a real convex
function on
$(\sum\limits ^n_{i=1} \oplus X_i)_p$ (it suffices to assume
that $f$ is defined on conv $\Omega$).	Let $\sigma _p$ be the Lipschitz
constant of $f$. Then, for all $c>0$,
$$P(|f-M_f| > c) \le 4 e^{-c^p/4\sigma ^p_p}\leqno (8)$$
\noindent where $M_f$ is the median of $f$. A similar inequality (with absolute
constants replacing the two fours) holds with expectation replacing the median:
$$P(|f-\bE f| > c) \le K e^{-\delta c^p/\sigma ^p_p}.\leqno (9)$$
\noindent (One can take $K=8, \delta = {1\over 32})$.

The proof of the first assertion is identical to the proof of Theorem 3 in
Talagrand's paper [T]. The second assertion follows from the first; see [MS], p.
142.

\noindent {\bf Remark 5:} Inequality (9) easily extends to the more general
setting where each $P_i$ is a Radon probability on $B_{X_i}$.

\noindent {\bf Remark 6:} In inequalities (8) and (9) the left hand side
involves
only the values of $f$ on $\Omega$ while the right hand side involves, through
$\sigma _p$, the values of $f$ on conv $\Omega$. Thus one can replace $\sigma
_p$ by the infinum of the Lipschitz constants of all convex extensions of
$f| _\Omega$ to conv $\Omega$.	We do not know how to compute this infimum even
in the original setup of Talagrand's theorem  where each $\Omega _i$ is a two
point set.

\noindent {An application.}

\proclaim Lemma 7. Let $\mu$ be a probability measure on $\{1,\ldots, N\}$. Let
$0 < r < s \le 2r$ and let $X$ be an $n$-dimensional subspace of
$L_r(\{1,\ldots, N\}, \mu)$ such that $\|x\|_s \le K\| x\|_r$ for all $x\in
X$ ($\|\cdot \|_s$ denotes the $L_s(\{1,\ldots, N\}, \mu)$ norm). Assume
moreover that $\mu(i) \le {2\over N}$ for all $1 \le i \le N$.	Then, for all
$1 > \varepsilon >0$  and all $k\ge c \varepsilon ^{-r} r^{1/p}(\log 
{2\over\varepsilon})^{1/p} K^r n^{1/p} N^{1/q}$, where 
$q = {s\over r},$ $p = {q \over q-1}$, there exists a subset 
$A \subseteq \{1,\ldots, N\}$ 
of cardinality k such that the restriction to $A$ is a multiple of a 
$(1+\varepsilon$)-isomorphism on $X$. In
particular, $X$  $(1+\varepsilon)$-embeds into $\ell^k_r$.

\noindent {\bf Proof:} Let $\delta _i, i=1,\ldots, N$, be independent mean
$\delta$
0,1-valued random variables. Fix \break $x \in X, \|x\|_r = 1$ and define
$f: \bR^N \to \bR$ by

$$f((a_i)^N_{i=1}) = \sum ^N_{i=1} a_i \mu(i)|x(i)|^r.$$

\noindent Then

$$\eqalign{\sigma _p(f) &= \sup _{\sum ^N_{i=1} |a_i|^p=1} \ \ \sum ^N_{i=1} a_i
\mu(i) |x(i)|^r\cr
&= \bigg( \sum ^N_{i=1} \mu(i)^q |x(i)|^s\bigg)^{1/q} \le \Big({2\over
N}\Big)^{q-1\over q} \bigg( \sum ^N_{i=1} \mu(i)|x(i)|^s\bigg) ^{r\over s}\cr
&\le \Big( {2\over N}\Big)^{1/p} K^r.}$$

\noindent It follows from Corollary 4 that

$$P\bigg( \bigg| \sum ^N_{i=1} \delta _i \mu(i) |x(i)|^r - \delta \bigg| >
c\bigg) \le 4 e^{-c^pN/8K^{rp}}\leqno (10)$$

\noindent and, consequently, using the usual estimate on the size of an
$\varepsilon$-net in $\partial B_X$; cf. [MS] p. 7, that

$$P\bigg( \bigg| \sum ^N_{i=1} \delta _i \mu(i)|x(i)|^r - \delta \bigg| \le
\varepsilon ^r\delta \ {\rm for \ all} \ x \in \partial B_X\bigg) \ge 1-4\exp
\Big(n r \log {2 \over \varepsilon} - \varepsilon ^{rp} \delta ^pN/8
K^{rp}\Big).$$

\noindent Set $k = 2 \delta N$ (= twice the average size of $\{ i; \delta _i =
1\}$). Then, for $\eta = c \varepsilon ^{rp} \Big(r \log {2\over
\varepsilon}\Big)^{-1}$ ($c$ universal) and $n \le \eta \delta ^p N/K^{rp}$, the
probability above is larger then ${1\over 2}$, so we can find a set of
cardinality $k$ which satisfies the requirement. Eliminating $\delta$ from the
two equations $k = 2\delta N$ and $n = \eta \delta ^p N/K^{rp}$ we get

$$k \approx 2\eta ^{-1/p} n^{1/p} N^{1/q} K^r.$$

\proclaim Theorem 8. [BLM]: \ Let $0 < r < t \le 2$ and let $T, \varepsilon >0$.
Then there exists a constant $C = C(\varepsilon, T, r, t)$ such that any $n$
dimensional subspace $X$ of $L_r$ with type $t$ constant $K,
(1+\varepsilon)$-embeds into $\ell ^N_r$ with $N \le C \cdot n$.\medskip

\noindent {\bf Proof:} \ We may assume that $X \subseteq L^M_r$ for some finite
$M$. By the Maurey-Nikishin-Rosenthal factorization theorem ([M] Theorem 8 
and Proposition 44), there exists a probability measure
$\mu$ on $\{1,\ldots, M\}$ such that $\|x\|_s \le K\|x\|_r$ for all $x\in X$
where $s = {r+t\over 2}$ and $K$ depends on $r, t$ and $T$ only. Splitting the
large atoms of $\mu$ into ones with measure $\le {1\over M}$ we get a new
probability measure $\bar \mu$ on $\{1,\ldots, N\}$ with $N \le 2M$ and $\bar
\mu(i) \le {1\over M} \le {2\over N}$.	$L_u(\mu)$ embeds isometrically (as a
sublattice) into $L_u(\bar \mu)$ for all $u$ and the inequality $\|x\|_s \le
K\|x\|_r$ for all $x\in X$ (in the new embedding) stays true. Applying Lemma 7
we get that $X \ (1+\delta)$-embeds into $L^k_r$ where $k$ is of order
$n^{1/p}M^{1/q}$.

Using the Maurey-Nikishin-Rosenthal theorem again we may repeat the argument 
to get that $X \ (1+\delta)^2$-embeds into $L^k_r$ for $k$ of order 
$n^{^{{1\over p}+{1\over pq}}}M^{^{1\over q^2}}$. Iterating one gets the 
result. This part is the same as in [BLM]. (One should be more careful than 
we have been above, taking the exact
form of $k$ into account, but it works.)

\noindent {\bf Remark 9:} Both B. Maurey and M. Talagrand pointed out to us
that versions of inequality (10) follow from known inequalities; in particular, 
(10) is an immediate consequence of the Azuma-Pisier inequality (see p. 45 in
[MS]) except that
the exponent on the right side of (10) must be  multiplied by a constant
$\delta_p$ which tends to infinity with p.  Since the degeneracy of 
this constant is unimportant for
proving Theorem 8, we in fact do not have a good application 
of (our slight generalization of) Talagrand's isoperimetric inequality. 
On the other hand, it is possible that the approach outlined above can be
used for general subspaces of $L_r$, in which case one expects to use a 
version of Lemma 7 with "s" close to "r", which forces "p" to infinity.

\centerline {\bf References}

\item{[BLM]} J. Bourgain, J. Lindenstrauss and V. Milman, Approximation of
zonoids by zonotopes, Acta Math., to appear.

\item{[M]} B. Maurey, Th\'eor\`emes de Factorisation pour les Op\'erateurs \`a
Valeurs dans un Espace $L^p$, Ast\'erisque Soc. Math. France, n$^\circ$ 11
(1974).

\item{[MS]} V. Milman and G. Schechtman, Asymptotic Theory of Finite Dimensional
Normed Spaces, Lecture Notes in Math. Vol 1200, Springer (1986).

\item{[T]} M. Talagrand, An isoperimetric theorem on the cube and the
Khinchine-Kahane inequalities, Proc. AMS (to appear).

\end